\documentclass[11pt,a4paper]{amsart}

\usepackage[utf8]{inputenc}
\usepackage{latexsym}
\usepackage{amssymb}
\usepackage{amsmath}

\title[Frequent hypercyclicity of random entire functions]{Frequent hypercyclicity of random entire functions for the differentiation operator}
\author[M. Nikula]{Miika Nikula}
\address{University of Helsinki, Department of Mathematics and Statistics,
         P.O. Box 68 , FIN-00014 University of Helsinki, Finland}
\email{miika.nikula@helsinki.fi}
\date{\today}

\keywords{Frequently hypercyclic operator, differentiation operator, rate of growth, entire functions, random construction}

\newtheorem{parg}{§}

\newtheorem{thm}[parg]{Theorem}
\newtheorem{prop}[parg]{Proposition}
\newtheorem{lemma}[parg]{Lemma}

\newtheorem{thmletter}{Theorem}

\newcommand{\N}{\mathbb{N}}
\newcommand{\Nz}{\mathbb{N}_0}

\newcommand{\C}{\mathbb{C}}
\newcommand{\E}{\mathbb{E}}

\newcommand{\Prob}{\mathbb{P}}

\newcommand{\ind}{\mathbf{1}}

\renewcommand{\d}{\, \mathrm{d}}
\renewcommand{\Re}{\mathrm{Re} \,}
\renewcommand{\Im}{\mathrm{Im} \,}

\newcommand{\eqlaw}{\stackrel{d}{=}}

\begin{document}

\maketitle

\begin{abstract}
In this note we study the random entire functions defined as power series $f(z) = \sum_{n=0}^\infty \frac{X_n}{n!} z^n$ with independent and identically distributed coefficients $(X_n)$ and show that, under very weak assumptions, they are frequently hypercyclic for the differentiation operator $D: H(\C) \to H(\C)$, $f \mapsto Df = f'$. This gives a very simple probabilistic construction of $D$-frequently hypercyclic functions in $H(\C)$. Moreover we show that, under more restrictive assumptions on the distribution of the $(X_n)$, these random entire functions have a growth rate that differs from the slowest growth rate possible for $D$-frequently hypercyclic entire functions at most by a factor of a power of a logarithm.
\end{abstract}

\section{Introduction}

\subsection{Definitions and Background}

Let $V$ be a topological vector space and $T: V \to V$ a linear operator. We call the operator $T$ \emph{hypercyclic} if there exists a $f \in V$ such that $\{T^n f: n \in \N\}$, the \emph{orbit} of $f$ under $T$, is dense in $V$. In this case we also say that $f$ is hypercyclic for $T$.

The stronger notion of \emph{frequent hypercyclicity} is defined as follows. The vector $f$ is frequently hypercyclic for $T$, if for any open set $U \subset V$ the sequence of iterates of $f$ under $T$ that belong to $U$ has a positive lower density, or explicitly
$$
\liminf_{n \to \infty} \frac{ \left\{ k \in \{0,1,\dots,n-1\} \,\big\vert\, T^k f \in U \right\} }{n} > 0.
$$

The topological vector space we study in this note is the space of entire functions
$$
H(\C) = \left\{ f: \C \to \C \,\big\vert\, f \textrm{ is holomorphic everywhere} \right\},
$$
equipped with the standard topology given by
$$
f_n \underset{n\to\infty}{\longrightarrow} f \textrm{ in } H(\C) \quad
\Longleftrightarrow
\quad f_n \underset{n\to\infty}{\longrightarrow} f \textrm{ uniformly on all compact sets } K \subset \C,
$$
and the operator is the differentiation operator
$$
D: H(\C) \to H(\C), \quad f \mapsto Df = f'.
$$
It was shown already by MacLane \cite{ma52} that, in this setting, the differentiation operator is hypercyclic. The concept of frequent hypercyclicity was defined rather recently by Bayart and Grivaux (see e.g. \cite{bagr04} and \cite{bagr06}), and the frequent hypercyclicity of differentiation in the space of entire functions was proven soon after that by Bonilla and Grosse-Erdmann \cite{boge06}. We refer the reader to \cite{blboge10} for more background. In \cite{boge07}, Bonilla and Grosse-Erdmann posed the problem which motivates the probabilistic construction of this note: what is the slowest growth rate possible for a $D$-frequently hypercyclic entire function? As the final answer to this question, the following theorem was established by Drasin and Saksman \cite{drsa11} by a explicit constructions.

\begin{thmletter}\label{thm:optimal-growth}
For $f \in H(\C)$ and $r \geq 0$, define the circle maximum
$$
M_{f,\infty}(r) = \sup_{\theta \in [0,2\pi)} |f(re^{i\theta})|
$$
to measure the growth rate of $f$. Then for any $c > 0$ there exists a $D$-frequently hypercyclic $f$ such that
$$
M_{f,\infty}(r) \leq c r^{-\frac{1}{4}} e^r \quad \textrm{for large } r.
$$
Conversely, any $D$-frequently hypercyclic $f \in H(\C)$ satisfies
$$
\limsup_{r \to \infty} r^\frac{1}{4} e^{-r} M_{f,p}(r) > 0.
$$
\end{thmletter}

Continuing the earlier work of Blasco, Bonilla and Grosse-Erdmann \cite{boboge10}, it was also proven in \cite{drsa11} that the optimal growth rates of the circle averages
$$
M_{f,p}(r) = \left( \int_0^{2\pi} |f(re^{i\theta})|^p \frac{\d \theta}{2 \pi} \right)^{1/p}, \quad p \geq 1.
$$
are given by
$$
M_{f,p}(r) \leq \begin{cases}
\varphi(r) r^{-\frac{1}{2}} e^r & \textrm{for } p = 1 \\
c r^{-\frac{1}{2p}} e^r & \textrm{for } 1 < p \leq 2 \\
c r^{-\frac{1}{4}} e^r & \textrm{for } 2 \leq p
\end{cases},
$$
where $\varphi:[0,\infty) \to [0,\infty)$ is an arbitrary nondecreasing function such that $\varphi(r) \to \infty$ as $r \to \infty$ and $c > 0$ is an arbitrary constant. The optimal result for $p=1$ was already obtained by Bonet and Bonilla in \cite{bobo12}. That these bounds are optimal means that the inequalities will fail for $D$-frequently hypercyclic functions if $\varphi$ is replaced by a constant or if $c$ is replaced by a nonincreasing function $\psi:[0,\infty) \to [0,\infty)$ such that $\psi(r) \to 0$ as $r \to \infty$.

The purpose of this note is to give a simple probabilistic construction of $D$-frequently hypercyclic functions in $H(\C)$ that grow almost as slowly as Theorem \ref{thm:optima-growth} permits. The construction and the $D$-frequent hypercyclicity of the resulting functions is stated as Theorem \ref{thm:hypercyclicity} and an estimate for their growth rate as Theorem \ref{thm:hypercyclic-growth}.

\begin{thm}\label{thm:hypercyclicity}
Let $X$ be a complex random variable such that the support of the law of $X$ is whole $\C$ and such that the decay condition
\begin{equation}\label{eq:hypercyclicity-tail-rate}
\textrm{for some } \beta > 0, \quad \limsup_{r \to \infty} \, (\log r)^{1+\beta} \Prob(|X| \geq r) < \infty
\end{equation}
is satisfied. Let $(X_n)_{n=0}^\infty$ be a sequence of independent random variables with the law of $X$ and denote
$$
f(z) = \sum_{n=0}^\infty \frac{X_n}{n!} z^n.
$$
Then for almost all realizations of the sequence $(X_n)_{n=0}^\infty$, the power series above represents an entire function which is frequently hypercyclic for the differentiation operator in $H(\C)$.
\end{thm}

Our growth estimate referred to above does not use (frequent) hypercyclicity in any way and is a consequence of the following general estimate for the growth of a random entire function.

\begin{prop}\label{prop:growth}
Let $X$ be a complex random variable such that the decay condition
\begin{equation}\label{eq:subgaussian-tail-rate}
\textrm{for some } C > 0, \quad \limsup_{t \to \infty} e^{C t^2} \, \E e^{t |X|} < \infty
\end{equation}
is satisfied. Let $(X_n)_{n=0}^\infty$ be a sequence of independent random variables with the law of $X$ and denote
$$
f(z) = \sum_{n=0}^\infty \frac{X_n}{n!} z^n.
$$
There exists a deterministic constant $C > 0$ such that almost surely the function $f$ satisfies
\begin{equation}\label{eq:growth}
\sup_{|z|=r} \left| f(z) \right| \leq C \sqrt{\log r} \frac{e^r}{r^{1/4}}
\end{equation}
for all sufficiently large $r$.
\end{prop}

Combining Theorem \ref{thm:hypercyclicity} with Proposition \ref{prop:growth} gives the main result of this note.

\begin{thm}\label{thm:hypercyclic-growth}
Let $X$ be a complex random variable such that the support of the law of $X$ is whole $\C$ and that the decay condition \eqref{eq:subgaussian-tail-rate} is satisfied. Then the random power series
$$
f(z) = \sum_{n=0}^\infty \frac{X_n}{n!} z^n
$$
almost surely represents a $D$-frequently hypercyclic entire function that satisfies the growth estimate \eqref{eq:growth}.
\end{thm}

\vspace{0.5cm}

\noindent To finish this introductory section we sketch the main idea of the argument for why the random power series in Theorem \ref{thm:hypercyclicity} represent $D$-frequently hypercyclic entire functions. The details of the proof and the study of the growth rate are presented in the following two sections.

\subsection{Sketch of the Main Argument}

Let $X$ be a complex random variable such that the support of its law is the whole plane\footnote{This requirement is clearly necessary for the power series \eqref{eq:powerseries} to represent a $D$-frequently hypercyclic entire function.}. Consider an infinite sequence $(X_n)_{n=0}^\infty$ of independent copies of $X$ and define the random function $f: \C \to \C$ by
\begin{equation}\label{eq:powerseries}
f(z) = \sum_{n=0}^\infty X_n \frac{z^n}{n!}.
\end{equation}
Under the condition
$$
\limsup_{n\to\infty} \left(\frac{X_n}{n!}\right)^\frac{1}{n} = 0 \quad \textrm{almost surely},
$$
which is a rather weak assumption on the tail of $|X|$, the power series \eqref{eq:powerseries} defines a random entire function. Under the additional but still weak assumptions on the tail of $|X|$ in Theorem \ref{thm:hypercyclicity}, this random entire function can be shown to be frequently hypercyclic for the differentiation operator.

The complex plane can be covered by a countable number of disks centered at the origin, and in any such disk the set of polynomials with rational coefficients is a countable dense set in the space of bounded analytic functions of the disk. To show the frequent hypercyclicity of our $f$, it is thus sufficient to show that for any disk $B_r$ of radius $r>0$ centered at the origin, any polynomial $p: \C \to \C$ and any $\varepsilon > 0$, almost surely we have
$$
\liminf_{n \to \infty} \frac{ \left\{ k \in \{0,1,\dots,n-1\} \,\big\vert\, \sup_{z \in B_r} |f^{(k)}(z)-p(z)| < \varepsilon \right\} }{n} > 0.
$$
To understand why this should be true, first note that as the event $\{|X-a|<\varepsilon\}$ has a positive probability for any $a \in \C$ and $\varepsilon > 0$ we have
$$
\Prob\left( \sup_{z \in B_r} |f(z)-p(z)| < \varepsilon \right) > 0.
$$
But by differentiation we have
$$
f^{(k)}(z) = \sum_{n=0}^\infty X_{n+k} \frac{z^n}{n!} \quad \textrm{for all } k \in \N,
$$
which implies the equality in distribution
$$
f^{(k)} \eqlaw f.
$$
Thus also
$$
\Prob\left( \sup_{z \in B_r} |f^{(k)}(z)-p(z)| < \varepsilon \right) > 0 \quad \textrm{for all } k \in \N.
$$
Moreover the events
$$
\left\{\sup_{z \in B_r} |f^{(k)}(z)-p(z)| < \varepsilon\right\} \quad \textrm{and} \quad \left\{\sup_{z \in B_r} |f^{(j)}(z)-p(z)| < \varepsilon\right\}
$$
can be expected to be, in an appropriate sense, almost independent if the difference of $k$ and $j$ is large, since up to a small perturbation the values of an analytic function in a disk of radius $r$ are determined by some finite number of power series coefficients. To conclude, it remains to observe that in an infinite sequence of independent trials an event will occur with frequency equal to the probability of the event.

Thus for the proof of Theorem \ref{thm:hypercyclicity} the main challenge is to make quantitative the sense in which the events referred to above are "almost" independent. It is in this quantitative analysis that the assumptions on the tail of $|X|$ enter the consideration.

\section{Frequent Hypercyclicity of Random Entire Functions}

Before embarking on the proof of Theorem \ref{thm:hypercyclicity} we present a simple lemma that shows that the assumed decay condition cannot be weakened much without losing not only the $D$-frequent hypercyclicity but even the infinite radius of convergence of the power series.

\begin{lemma}\label{lemma:entire-power-series}
Let $X$ be a complex random variable such that
\begin{equation}\label{eq:divergence-tail-rate}
\liminf_{r \to \infty} \, (\log r) \Prob(|X| \geq r) > 0
\end{equation}
and let $(X_n)_{n=0}^\infty$ be a sequence of independent copies of $X$. Then almost surely the radius of convergence of the power series 
\begin{equation}\label{eq:powerseries-lemma}
\sum_{n=0}^\infty \frac{X_n}{n!} z^n
\end{equation}
is zero.
\end{lemma}

\emph{Remark.} The condition \eqref{eq:divergence-tail-rate} cannot be weakened much. For example, an easy argument similar to the proof of the lemma shows that if
$$
\limsup_{r \to \infty} \, (\log r) (\log \log r)^\alpha \Prob(|X| \geq r) < \infty \quad \textrm{for some } \alpha > 0
$$
the lemma fails and the power series \eqref{eq:powerseries-lemma} almost surely represents an entire function.

\vspace{0.5cm}

\emph{Proof.} Let $M > 0$. From \eqref{eq:divergence-tail-rate} we get, for some $c > 0$ and large enough $n$,
$$
\Prob\left(|X_n| \geq (M n)^n\right) \geq \frac{c}{n \log (M n)}.
$$
The series $\sum_{n=2}^\infty \frac{1}{n \log n}$ is divergent and the random variables $(X_n)_{n=1}^\infty$ are independent, so by the Borel--Cantelli lemma, almost surely there exists an infinite increasing sequence $(n_k)_{k=1}^\infty$ of indices such that
$$
\textrm{for all }k, \quad |X_{n_k}| \geq (M n_k)^{n_k}.
$$
The radius of convergence $R$ of the power series \eqref{eq:powerseries-lemma} is given by
$$
\frac{1}{R} = \limsup_{n \to \infty} \left(\frac{|X_n|}{n!}\right)^\frac{1}{n}.
$$
By using Stirling's approximation for $n!$, we get
$$
\left(\frac{|X_{n_k}|}{n_k!}\right)^\frac{1}{n_k} \sim \frac{e}{n_k} |X_{n_k}|^\frac{1}{n_k} (2\pi n_k)^{-\frac{1}{2n_k}} \geq e M (2 \pi n_k)^{-\frac{1}{2n_k}}
$$
for the sequence $(n_k)_{k=1}^\infty$. It follows that
$$
\frac{1}{R} \geq e M \quad \textrm{a.s.}
$$
that is, almost surely the radius of convergence of \eqref{eq:powerseries-lemma} is less than $\frac{1}{eM}$.

Since $M > 0$ was arbitrary, it follows that the radius of convergence of \eqref{eq:powerseries-lemma} is almost surely zero. $\Box$

\vspace{0.5cm}

By the remark after the statement of Lemma \ref{lemma:entire-power-series} we now know that the condition \eqref{eq:hypercyclicity-tail-rate} ensures that the random power series $\sum_{n=0}^\infty \frac{X_n}{n!} z^n$ almost surely defines an entire function. We are thus ready to start the proof of its hypercyclicity.

\vspace{0.5cm}

\emph{Proof of Theorem \ref{thm:hypercyclicity}.} Let $r > 0$ and $\varepsilon > 0$ be fixed and let
$$
g(z) = \sum_{n=0}^N a_n \frac{z^n}{n!}
$$
be a given polynomial of degree $N$. As indicated in the introduction, to prove the theorem it is sufficient to show that almost surely,
$$
\sup_{z \in B_r} \left| f^{(k)}(z) - g(z) \right| < \varepsilon
$$
for a set of indices $k$ with a positive lower density. Defining the events
$$
A_k = \left\{ \sup_{z \in B_r} \left| f^{(k)}(z) - g(z) \right| < \varepsilon \right\}
$$
and the random variable
$$
S_n' = \sum_{k=0}^{n-1} \ind_{A_k},
$$
this is equivalent to showing that
$$
\liminf_{n \to \infty} \frac{S_n'}{n} > 0.
$$
A direct approach to this kind of a problem is to compute the expectation and variance of $S_n'$. If the expectation grows linearly but variance slower than quadratically, we are done.

Computation of the variance of $S_n'$ involves the conditional probabilities $\Prob(A_j|A_k)$ for $j,k \in \Nz$, which are rather difficult to get to directly.
To get more manageable quantities we define the events $B_k$ by setting, for all $k \in \Nz$,
$$
B_k = \left\{ \sum_{n=0}^N |X_{k+n}-a_n| \frac{r^n}{n!} < \frac{\varepsilon}{2} \quad \textrm{and} \quad |X_{k+n}| < \rho_n \textrm{ for all }n > N \right\}
$$
where $(\rho_n)_{n=N+1}^\infty$ is a nondecreasing sequence of positive reals tending to infinity such that
\begin{equation}\label{eq:rho-upper-growth-condition}
\sum_{n=N+1}^\infty \rho_n \frac{r^n}{n!} < \frac{\varepsilon}{2}
\end{equation}
and that
\begin{equation}\label{eq:rho-lower-growth-condition}
\sum_{n=d}^\infty \Prob(|X| \geq \rho_n) = O(d^{-\gamma}) \quad \textrm{for some } \gamma > 0.
\end{equation}
We postpone the choice of the sequence $(\rho_n)_{n=N+1}^\infty$ until the end of the proof and first show how its existence implies the frequent hypercyclicity of $f$.

By \eqref{eq:rho-upper-growth-condition} and the definition of the events $B_k$ it is clear that $B_k \subset A_k$ for all $k \in \Nz$.
Thus
$$
S_n := \sum_{k=0}^{n-1} \ind_{B_k} \leq \sum_{k=0}^{n-1} \ind_{A_k} = S_n'.
$$
We will show that almost surely
\begin{equation}\label{eq:modified-event-positive-lower-density}
\liminf_{n \to \infty} \frac{S_n}{n} > 0,
\end{equation}
which implies the theorem.

Denote
$$
p = \Prob\left( \sum_{k=0}^N |X_k-a_k| \frac{r^k}{k!} < \frac{\varepsilon}{2} \right)
$$
and
$$
Q_d = \prod_{k=d}^\infty \Prob(|X| < \rho_j) = \prod_{k=d}^\infty \big(1-\Prob(|X| \geq \rho_j)\big) \quad \textrm{for }d = N+1,N+2,\dots.
$$
Since the $X_n$ are independent and the support of their law is the whole plane, we have $p > 0$. The assumption \eqref{eq:rho-lower-growth-condition} in turn implies that all the products $Q_d$ are also positive. Using these notations, the invariance of the law of $f$ under differentiation gives
\begin{eqnarray*}
\E S_n \, = \, \E \sum_{k=0}^{n-1} \ind_{B_k} & = & \sum_{k=0}^{n-1} \Prob(B_k) \, = \, n \Prob(B_0) \\
& = & n \, \Prob\left( \sum_{n=0}^N |X_n-a_n| \frac{r^n}{n!} < \frac{\varepsilon}{2} \right) \prod_{n=N+1}^\infty \Prob\left( |X_n| < \rho_n \right) \\
& = & n p Q_{N+1}.
\end{eqnarray*}
Similarly, the variance of $S_n$ is given by
\begin{eqnarray}
\E (S_n - \E S_n)^2 & = & \E \left(\sum_{k=0}^{n-1} \ind_{B_k}\right)^2 - (\E S_n)^2 \nonumber \\
& = & \sum_{k=0}^{n-1} \sum_{j=0}^{n-1} \Prob(B_k \cap B_j) - n^2 p^2 Q_{N+1}^2 \nonumber \\
& = & n p Q_{N+1} + 2 \sum_{d=1}^{n-1} (n-d) \Prob(B_0 \cap B_d) - n^2 p^2 Q_{N+1}^2. \label{eq:sn-variance-probabilities}
\end{eqnarray}
The probability $\Prob(B_0 \cap B_d)$ can be explicitly computed for large $d$. Since the sequence $(\rho_k)_{k=N+1}^\infty$ is nondecreasing and tends to infinity, there exists a constant $M \in \N$, $M > N$ the choice of which depends on $g$, $r$ and $(\rho_k)_{k=N+1}^\infty$ such that for all $d \geq M$,
\begin{equation}\label{eq:d-condition}
|a_k| + \frac{\varepsilon}{2} \frac{k!}{r^k} \leq \rho_{k+d} \quad \textrm{for all } k=0,1,\dots,N.
\end{equation}
This condition ensures that
$$
\left\{ \sum_{k=0}^N |X_{d+k}-a_k| \frac{r^k}{k!} < \frac{\varepsilon}{2} \right\} \subset \left\{ |X_{k+d}| \leq \rho_{k+d} \quad \textrm{for } k=0,1,\dots,N \right\}
$$
for all $d \geq M$. Thus by the definition of the events $B_j$, for all $d \geq M$ we have
\begin{eqnarray*}
\Prob(B_0 \cap B_d) & = & \Prob\left( \sum_{k=0}^N |X_k-a_k| \frac{r^k}{k!} < \frac{\varepsilon}{2}, \quad \sum_{k=0}^N |X_{d+k}-a_k| \frac{r^k}{k!} < \frac{\varepsilon}{2},  \right. \\
& & \quad |X_k| < \rho_k \textrm{ for }k > N \quad \textrm{and} \quad |X_{d+k}| < \rho_k \textrm{ for }k > N \Bigg) \\
& = & \Prob\left( \sum_{k=0}^N |X_k-a_k| \frac{r^k}{k!} < \frac{\varepsilon}{2},  \quad \sum_{k=0}^N |X_{d+k}-a_k| \frac{r^k}{k!} < \frac{\varepsilon}{2},  \right. \\
& & \quad |X_k| < \rho_k \textrm{ for } N < k < d \quad \textrm{and} \quad |X_{d+k}| < \rho_k \textrm{ for } k > N \Bigg) \\
& = & p^2 \prod_{k=N+1}^{d-1} \Prob(|X|<\rho_k) \prod_{k=N+1}^\infty \Prob(|X|<\rho_k) \, = \, p^2 \frac{Q_{N+1}^2}{Q_d}.
\end{eqnarray*}
Continuing from \eqref{eq:sn-variance-probabilities},
\begin{eqnarray*}
\E (S_n - \E S_n)^2 & = & n p Q_{N+1} + 2 \sum_{d=1}^{M-1} (n-d) \Prob(B_0 \cap B_d) \\
& & + \, 2 \sum_{d=M}^{n-1} (n-d) p^2 \frac{Q_{N+1}^2}{Q_d} - n^2 p^2 Q_{N+1}^2.
\end{eqnarray*}
for all $n > M$. We write the last term as
\begin{eqnarray*}
n^2 p^2 Q_{N+1}^2 & = & p^2 Q_{N+1}^2 \left( 2 \sum_{d=1}^{n-1} (n-d) + n \right) \\
& = & n p^2 Q_{N+1}^2 + 2 \left( \sum_{d=1}^{M-1} + \sum_{d=M}^{n-1} \right) (n-d) p^2 Q_{N+1}^2
\end{eqnarray*}
and regroup the terms to get
\begin{eqnarray}
\E (S_n - \E S_n)^2 & = & n (p Q_{N+1} - p^2 Q_{N+1}^2) + 2 \sum_{d=1}^{M-1} (n-d) \left( \Prob(B_0 \cap B_d) - p^2 Q_{N+1}^2 \right) \nonumber \\
& & + 2 \sum_{d=M}^{n-1} (n-d) \left( p^2 \frac{Q_{N+1}^2}{Q_d} - p^2 Q_{N+1}^2 \right). \label{eq:sn-variance}
\end{eqnarray}
Obviously the terms on the first line grow as $O(n)$ as $n$ tends to infinity. To get a bound on the asymptotic growth rate of the sum on the second line we use the elementary inequalities
$$
\frac{1}{1-x} \leq e^{c_1 x} \quad \textrm{and} \quad e^{c_1 x} - 1 \leq c_1 c_2 x
$$
that hold for some constants $c_1,c_2 > 0$ in some neighbourhood of the origin. By taking the constants large enough, these inequalities can be assumed to hold on the interval $[0,1/2]$. These inequalities allow us to estimate
\begin{eqnarray*}
\frac{1}{Q_d} - 1 & = & \prod_{k=d}^\infty \big(1-\Prob(|X| \geq \rho_k)\big)^{-1} - 1 \, \leq \, \prod_{k=d}^\infty e^{c_1 \Prob(|X| \geq \rho_k)} - 1 \\
& = & e^{c_1 \sum_{k=d}^\infty P(\rho_k)} - 1 \, \leq \, c_1 c_2 \sum_{k=d}^\infty \Prob(|X| \geq \rho_k)
\end{eqnarray*}
for $d$ large enough that $\sum_{k=d}^\infty P(\rho_k) \leq 1/2$. Recalling the choice of the constant $M$, we may suppose that $M$ has been taken large enough that this estimate holds for all $d \geq M$. We use this estimate in \eqref{eq:sn-variance} to get
$$
\E (S_n - \E S_n)^2 \leq O(n) + 2 p^2 Q_{N+1}^2 c_1 c_2 \sum_{d=M}^{n-1} (n-d) \sum_{k=d}^\infty \Prob(|X| \geq \rho_k),
$$
and by the property \eqref{eq:rho-lower-growth-condition} of the sequence $(\rho_k)_{k=N+1}^\infty$ this further implies
$$
\E (S_n - \E S_n)^2 \leq O(n) + O\left( \sum_{d=M}^{n-1} (n-d) \, d^{-\gamma} \right) \quad \textrm{for some } \gamma > 0.
$$
For $0 < \gamma < 1$ we determine the asymptotics of the sum by writing
$$
\sum_{d=M}^{n-1} (n-d) \, d^{-\gamma} = n^{2-\gamma} \sum_{d=M}^{n-1} \left( 1 - \frac{d}{n} \right) \left( \frac{d}{n} \right)^{-\gamma} \frac{1}{n}
$$
and noting that the last sum tends to the integral $\int_0^1 (1-x) x^{-\gamma} \d x < \infty$ as $n \to \infty$. For $0 < \gamma < 1$ we thus have
$$
\E (S_n - \E S_n)^2 = O(n^{2-\gamma}).
$$
For $\gamma \geq 1$ the asymptotic growth of the variance is even slower, since we may estimate
$$
\sum_{d=M}^{n-1} (n-d) \, d^{-\gamma} \leq n \sum_{d=M}^{n-1} \left(\frac{1}{d} - \frac{1}{n}\right) = O( n \log n ).
$$

Up to now we have deduced that under the assumptions \eqref{eq:rho-upper-growth-condition} and \eqref{eq:rho-lower-growth-condition} we have
$$
\E S_n = n p Q_{N+1} \quad \textrm{and} \quad \E (S_n - \E S_n)^2 = O(n^{2-\alpha}) \quad \textrm{for some } 0 < \alpha < 1.
$$
These facts imply that \eqref{eq:modified-event-positive-lower-density} holds almost surely. To see this we start by using Chebychev's inequality to get
$$
\Prob\left( \left|S_n - \E S_n\right| \geq n^{1-\frac{\alpha}{4}} \right) \leq n^{-2 \left(1-\frac{\alpha}{4}\right)} \E (S_n - \E S_n)^2 = O\left(n^{-\frac{\alpha}{2}}\right).
$$
Choose a $\delta > 1$ so that $-\frac{\delta \alpha}{2} < -1$ and consider a sequence $(n_j)_{j=1}^\infty$ of indices such that $n_j \sim j^\delta$. For this sequence the series
$$
\sum_{j=1}^\infty \Prob\left( \left|S_{n_j} - \E S_{n_j}\right| \geq n_j^{1-\frac{\alpha}{4}} \right)
$$
is summable, so by the Borel--Cantelli lemma, almost surely
$$
n_j^{1-\frac{\alpha}{4}} > |S_{n_j}-\E S_{n_j}| \quad \Longrightarrow \quad p Q_{N+1} - n_j^{-\frac{\alpha}{4}} < \frac{S_{n_j}}{n_j} \quad \textrm{for all but finitely many } j.
$$
This already implies that $\liminf_{j \to \infty} \frac{S_{n_j}}{n_j} > 0$ almost surely. To extend this from the subsequence $(n_j)$ to the full result, let $n \geq n_1$ be arbitrary and take an index $j$ so that $n_j \leq n \leq n_{j+1}$. Since $S_n$ is nondecreasing we have
$$
\frac{S_n}{n} > \frac{S_{n_j}}{n_{j+1}} = \frac{n_j}{n_{j+1}} \frac{S_{n_j}}{n_j} > \frac{n_j}{n_{j+1}} ( p Q_{N+1} - n_j^{-\frac{\alpha}{4}}).
$$
Noting that the asymptotic growth $n_j \sim j^\delta$ implies $\lim_{j \to \infty} \frac{n_j}{n_{j+1}} = 1$, the full claim \eqref{eq:modified-event-positive-lower-density} follows.

It remains to show that under the assumption \eqref{eq:hypercyclicity-tail-rate} it is possible to choose a sequence $(\rho_k)_{k=N+1}^\infty$ satisfying \eqref{eq:rho-upper-growth-condition} and \eqref{eq:rho-lower-growth-condition}. Our choice is $\rho_n = a e^n$, where $a > 0$ is taken so that
$$
a \sum_{n=N+1}^\infty \frac{(er)^n}{n!} < \frac{\varepsilon}{2}.
$$
This choice ensures that \eqref{eq:rho-upper-growth-condition} holds. To see that \eqref{eq:rho-lower-growth-condition} holds as well, observe that by \eqref{eq:hypercyclicity-tail-rate} there exists a constant $c > 0$ such that
$$
\Prob(|X| \geq a e^n) \leq \frac{c}{(n + \log a)^{1+\beta}}
$$
for all $n \in \N$ such that $a e^n \geq 1$. Thus for large $d$ we have
$$
\sum_{n=d}^\infty \Prob(|X| \geq a e^n) \leq \sum_{n=d}^\infty \frac{c}{(n + \log a)^{1+\beta}} \leq \int_{d-1+\log a}^\infty \frac{c}{x^{1+\beta}} \d x = O(d^{-\beta}),
$$
which means that \eqref{eq:rho-lower-growth-condition} holds with $\gamma = \beta$. The proof is complete. $\Box$

\vspace{0.5cm}

\emph{Remark.} In \cite{bama11} Bayart and Matheron give an eigenvalue characterization for an operator to be weakly mixing with respect to a Gaussian measure. If the distribution of $X$ in our construction is Gaussian, the measure we obtain on $H(\C)$ is clearly an explicit example of a measure of the type considered by Bayart and Matheron.

\section{Growth of Random Analytic Functions}

In this section we show that random frequently hypercyclic entire functions have slow growth rates among all possible frequently hypercyclic entire functions.

Let us first consider a random entire function $f$ defined as in Theorem \ref{thm:hypercyclicity} for a standard complex Gaussian random variable $X$. This means that $\Re X$ and $\Im X$ are independent real Gaussians with variance $\frac{1}{2}$. Since sums of independent Gaussian variables are Gaussian, the values $\{f(z)\}_{z \in \C}$ are also Gaussian random variables. Moreover, the variance of the sum of independent random variables is the sum of the variances, so we have
$$
\E |f(z)|^2 = \E \left| \sum_{n=0}^\infty \frac{X_n}{n!} z^n \right|^2 = \sum_{n=0}^\infty \frac{|z|^{2n}}{(n!)^2}.
$$

Let us denote
$$
I(r) = \sum_{n=0}^\infty \frac{r^{2n}}{(n!)^2} \quad \textrm{for } r \geq 0.
$$
The growth rates of our random analytic functions are determined by the asymptotics of the function $I(r)$. Our function $I(r)$ is, up to a scaling of the argument, the same as the modified Bessel function of the second kind, usually denoted by $I_0$; namely we have $I(r) = I_0(2r)$ for all $r \geq 0$. By the asymptotics of $I_0$ we have (see e.g. formula 9.7.1 in \cite{abst72})
$$
I(r) \sim \frac{1}{2 \sqrt{\pi}} \frac{e^{2r}}{\sqrt{r}} \quad \textrm{as } r \to \infty.
$$
The formula $\E |X|^p = \Gamma(\frac{p}{2}+1) (\E |X|^2)^{p/2}$ for the moments of a Gaussian random variable now gives, for all $p > 0$,
\begin{eqnarray*}
\E \int_0^{2 \pi} |f(r e^{i \theta})|^p \frac{\d \theta}{2 \pi} & = & \int_0^{2 \pi} \E |f(r e^{i \theta})|^p \frac{\d \theta}{2 \pi} \, = \, \E |f(r)|^p \\
& = & \Gamma\left(\frac{p}{2}+1\right) I(r)^{\frac{p}{2}} \\
& \asymp & \left( \frac{e^r}{r^{\frac{1}{4}}} \right)^p \quad \textrm{as } r \to \infty.
\end{eqnarray*}
Recalling Theorem \ref{thm:optimal-growth}, this calculation indicates that the integral means of the random analytic functions in Theorem \ref{thm:hypercyclicity} have a slow rate of growth among all $D$-frequently hypercyclic functions, at least when averaged over the possible realizations of the randomness. In Proposition \ref{prop:growth} we consider the individual realizations of the construction.

The proof of Proposition \ref{prop:growth} relies on the following estimate of Kahane (\cite{ka85}) for the probability of a random trigonometric polynomial taking large values.

\begin{lemma}\label{lemma:kahane-estimate}
Let $\xi$ be a real random variable that satisfies the decay condition \eqref{eq:subgaussian-tail-rate}. Let $(\xi_n)_{n=1}^K$ be a finite sequence of independent copies of $\xi$, $(a_n)_{n=1}^K$ a finite sequence of complex constants and $(q_n)_{n=1}^K$ a finite sequence of trigonometric polynomials of degree less than or equal to $N$. Then there exists a constant $c > 0$ that depends only on the distribution of $\xi$ such that the tail estimate
\begin{equation}\label{eq:kahane-estimate}
\Prob\left( \sup_{\theta \in [0,2\pi]} \left| \sum_{n=1}^K a_n \xi_n q_n(\theta) \right| > c \sqrt{\log N} \, \left( \sum_{n=1}^K |a_n|^2 \right)^\frac{1}{2} \right) \leq \frac{1}{N^2}
\end{equation}
holds.
\end{lemma}

\emph{Proof of Proposition \ref{prop:growth}.} Write
$$
f(z) \, = \, \sum_{n=0}^\infty X_n \frac{z^n}{n!} \, = \, \sum_{n=0}^\infty (\Re X_n) \frac{z^n}{n!} + i \sum_{n=0}^\infty (\Im X_n) \frac{z^n}{n!} =: g(z) + h(z).
$$
It is clearly enough to prove that the conclusion of the proposition holds for the function $g$.

Let $R \in \N$ and further decompose $g$ by writing
\begin{eqnarray*}
g(z) \, = \, \sum_{n=0}^\infty (\Re X_n) \frac{z^n}{n!} & = & \sum_{n=0}^{3R} (\Re X_n) \frac{z^n}{n!} + \sum_{j=1}^\infty \sum_{n=3^j R+1}^{3^{j+1} R} (\Re X_n) \frac{z^n}{n!} \\
& =: & g_{R,0}(z) + \sum_{j=1}^\infty g_{R,j}(z).
\end{eqnarray*}
We will show that almost surely for all but finitely many $R \in \N$, around $|z|=R$ the contribution of the series $\sum_{j=1}^\infty g_{R,j}$ is negligible and that $g_{R,0}$ has the claimed order of magnitude.

We start with the analysis of $g_{R,0}$. As
$$
\sup_{|z|=R} |g_{R,0}(z)| = \sup_{\theta \in [0,2\pi]} \left| \sum_{n=0}^{3R} X_n \frac{R^n}{n!} e^{in\theta} \right|
$$
and
$$
\sum_{n=0}^{3R} \left( \frac{R^n}{n!} \right)^2 \leq \sum_{n=0}^\infty \frac{R^{2n}}{(n!)^2} = I(R),
$$
by Lemma \ref{lemma:kahane-estimate} there exists a constant $c_1$ depending only on the distribution of $\Re X$ such that
$$
\Prob\left( \sup_{|z|=R} |g_{R,0}(z)| > c_1 \sqrt{\log 3R} \, \sqrt{I(R)} \right) \leq \frac{1}{(3R)^2}.
$$
The series $\sum_{R=1}^\infty \frac{1}{R^2}$ is summable, so by the Borel--Cantelli lemma we have
\begin{equation}\label{eq:principal-growth}
\sup_{|z|=R} |g_{R,0}(z)| \leq c_1 \sqrt{\log 3R} \, \sqrt{I(R)}
\end{equation}
almost surely for all but finitely many $R \in \N$. Since
$$
\sqrt{I(R)} \asymp \frac{e^R}{R^{\frac{1}{4}}}
$$
as $R \to \infty$, we see that indeed $g_{R,0}(z)$ has the claimed order of magnitude close to $|z| = R$.

For treating the functions $g_{R,j}$ we need an estimate for the tail of the series of $I(R)$. By bounding from above by a geometric series we see that for any integer $A \geq 3R$ we have
$$
\sum_{n=A+1}^\infty \frac{R^{2n}}{(n!)^2} < \frac{R^{2A}}{(A!)^2} \sum_{n=1}^\infty \frac{R^{2n}}{A^{2n}} = \frac{R^{2A}}{(A!)^2} \frac{(R/A)^2}{1-(R/A)^2} \leq \frac{1}{8} \frac{R^{2A}}{(A!)^2}.
$$
Applying this and Lemma \ref{lemma:kahane-estimate} to the functions $g_{R,j}$ gives
\begin{eqnarray*}
\Prob\left( \sup_{|z|=R} |g_{R,j}(z)| > c_1 \sqrt{\log 3^{j+1}R} \frac{1}{8} \frac{R^{2 \cdot 3^j R}}{((3^j R)!)^2} \right) \leq \frac{1}{(3^{j+1}R)^2}.
\end{eqnarray*}
The double series $\sum_{j=1}^\infty \sum_{R=1}^\infty \frac{1}{(3^{j+1}R)^2}$ is summable. The Borel--Cantelli lemma thus implies that almost surely
$$
\sup_{|z|=R} |g_{R,j}(z)| \leq c_1 \sqrt{\log 3^{j+1}R} \frac{1}{8} \frac{R^{2 \cdot 3^j R}}{((3^j R)!)^2}
$$
for all but finitely many pairs $(R,j) \in \N^2$. We use Stirling's approximation to write, for large $R$,
\begin{eqnarray*}
c_1 \sqrt{\log 3^{j+1}R} \frac{1}{8} \frac{R^{2 \cdot 3^j R}}{((3^j R)!)^2} & \sim & \frac{c_1}{8} \sqrt{(j+1)\log 3 + \log R} \frac{R^{2 \cdot 3^j R}}{\left(\frac{3^jR}{e}\right)^{2 \cdot 3^jR} 2 \pi \, 3^j R } \\
& = & \frac{c_1}{16 \pi} \sqrt{(j+1)\log 3 + \log R} \left( \frac{e}{3^j} \right)^{2 \cdot 3^j R} \frac{1}{3^j R}.
\end{eqnarray*}
From this form of the estimate for $\sup_{|z|=R} |g_{R,j}(z)|$ it is easy to see that we have actually shown that almost surely
\begin{equation}\label{eq:growth-negligible-part}
\sum_{j=1}^\infty \sup_{|z|=R} |g_{R,j}(z)| \to 0 \quad \textrm{as } R \to \infty.
\end{equation}

The estimates \eqref{eq:principal-growth} and \eqref{eq:growth-negligible-part} show that almost surely, for all but finitely many $R \in \N$ we have
$$
\sup_{|z|=R} |g(z)| \leq \sup_{|z|=R} |g_{R,0}(z)| + \sum_{j=1}^\infty \sup_{|z|=R} |g_{R,j}(z)| \leq C \sqrt{\log R} \frac{e^R}{R^\frac{1}{4}}
$$
for some deterministic constant $C>0$ that depends only on the distribution of $\Re X$. The extension from integer to real radii is immediate, as the ratio
$$
\frac{C \sqrt{\log (R+1)} \frac{e^{R+1}}{(R+1)^\frac{1}{4}}}{C \sqrt{\log R} \frac{e^R}{R^\frac{1}{4}}} = e \sqrt{ \frac{\log(R+1)}{\log R} } \left( \frac{R}{R+1} \right)^\frac{1}{4}
$$
is bounded as $R \to \infty$ and by the maximum principle, the function $R \mapsto \sup_{|z|=R} |g(z)|$ is increasing. $\Box$

\vspace{0.5cm}

\emph{Remark.} If the random variable $X$ is Gaussian, a somewhat more elaborate argument can be used to show that the growth estimate \eqref{eq:growth} is optimal.

\vspace{1.0cm}

\noindent \emph{Acknowledgements.} I thank David Drasin and my advisor Eero Saksman for helpful comments and for directing me to this problem.

\end{document}